.
.
.
\font\sets=msbm10.
\font\stampatello=cmcsc10.
\font\symbols=msam10.

\def\1{{\bf 1}}
\def\sgn{{\rm sgn}}

\def\dashsum{\mathop{\enspace{\sum}'}}
\def\starsum{\mathop{\enspace{\sum}^{\ast}}}

\def\square{\hbox{\vrule\vbox{\hrule\phantom{s}\hrule}\vrule}}
\def\defineq{\buildrel{def}\over{=}}
\def\defin{\buildrel{def}\over{\Longleftrightarrow}}
\def\doublesum{\mathop{\sum\sum}}
\def\doublestarsum{\mathop{\starsum\starsum}}

\def\supporto{{\rm supp}\,}

\def\N{\hbox{\sets N}}

\def\R{\hbox{\sets R}}
\def\Z{\hbox{\sets Z}}

\def\EssBdd{\hbox{\symbols n}\,}

\par
\centerline{\bf ON THE SYMMETRY INTEGRAL}
\bigskip
\centerline{\stampatello giovanni coppola}
\bigskip
\rightline{\it to my Angel}
\bigskip
{
\font\eightrm=cmr8
\eightrm {
\par
{\stampatello abstract.} We give a level one result for the \lq \lq symmetry integral\rq \rq, say $I_f(N,h)$, of essentially bounded $f:\N \rightarrow \R$; i.e., we get a kind of \lq \lq square-root cancellation\rq \rq \thinspace bound for the mean-square (in $N<x\le 2N$) of the \lq \lq symmetry sum\rq \rq \thinspace of, say, the arithmetic function $f:=g\ast \1$, where $g:\N \rightarrow \R$ is such that $\forall \varepsilon>0$ we have $g(n)\ll_{\varepsilon} n^{\varepsilon}$, and supported in $[1,Q]$, with $Q\ll N$ (so, the exponent of $Q$ relative to $N$, say the level $\lambda:=(\log Q)/(\log N)$ is $\lambda < 1$), where the symmetry sum weights the $f-$values in (almost all, i.e. all but $o(N)$ possible exceptions) the short intervals $[x-h,x+h]$ (with positive/negative sign at the right/left of $x$), with mild restrictions on $h$ (say, $h\to \infty$ and $h=o(\sqrt N)$, as $N\to \infty$).
}
\footnote{}{\par \noindent {\it Mathematics Subject Classification} $(2010) : 11N37, 11N25.$}
}
\bigskip
\par
\noindent \centerline{\stampatello 1. Introduction and statement of the results.}
\smallskip
\par
\noindent
We give upper bounds for the {\it symmetry integral} (compare [C-S], [C], [C1], [C2], [C3], [C5], [C6], esp.)
$$
I_f(N,h)\defineq \sum_{x\sim N} \Big| \dashsum_{|n-x|\le h} \sgn(n-x)f(n)\Big|^2, 
$$
\par
\noindent
with $x\sim N$ for $N<x\le 2N$, $\sgn(r):={r\over {|r|}}$ $\forall r\in \R$, $r\neq 0$, $\sgn(0):=0$ and the dash means : the terms $n=x\pm h$ have to be halved (say, $\dashsum_{|n-x|\le h} \sgn(n-x)f(n)=\sum_{|n-x|\le h}\sgn(n-x)f(n)-{{f(x+h)-f(x-h)}\over 2}$); this is one of the possible definitions : comparing the one(s) in the papers quoted above, the difference involved will be soon shown negligible, at least for the (large) class of {\it essentially bounded} arithmetic functions $f:\N \rightarrow \R$. Here we'll use the term \lq \lq {\it essentially}\rq \rq \thinspace to leave (negligible) multiplicative factors bounded by arbitrary small powers of $N$. For example, 
$$
f \enspace \hbox{\rm is } \hbox{\stampatello essentially bounded } (\hbox{\rm abbrev. } f\EssBdd 1) \enspace \defin \enspace \forall \varepsilon>0 \enspace f(n)\ll_{\varepsilon} N^{\varepsilon}
$$
\par
\noindent
(here, $\forall n\le 2N+h$, since we \lq \lq don't see\rq \rq \thinspace $f$ any further), with the usual Vinogradov notation (i.e., $\ll_{\varepsilon}$ means \lq \lq bounded in absolute value with a constant factor depending on $\varepsilon$\rq \rq); and we'll abbreviate 
$$
F(N,h)\EssBdd G(N,h) \enspace \defin \enspace \forall \varepsilon>0 \enspace |F(N,h)|\ll_{\varepsilon} N^{\varepsilon}G(N,h).
$$
\bigskip
\par
\noindent
For example, to see that our new definition is \lq \lq close\rq \rq \thinspace to the old one (an integral !), whenever $f\EssBdd 1$, 
$$
\int_{N}^{2N}\Big| \sum_{|n-x|\le h} \sgn(n-x)f(n)\Big|^2\,dx\EssBdd \sum_{N\le x<2N}\Big| \dashsum_{|n-[x]|\le h} \sgn(n-[x])f(n)\Big|^2\,dx + N 
\EssBdd I_f(N,h) + N + h^2 
$$
\par
\noindent
where $[x]\defineq$ {\it integer part} of $x\in \R$, both the remainders $\EssBdd N$ and $\EssBdd h^2$ are clearly negligible : below the \lq \lq {\it diagonal}\rq \rq \thinspace remainders, i.e. $\EssBdd Nh$ (when $h\to \infty$ and $h=o(N)$, for $N\to \infty$, as we'll assume henceforth). 
\par
The symmetry integral (now on, $I_f(N,h)$ defined above) measures, for $f$, the {\stampatello almost-all} (i.e., leave eventual $o(N)$ exceptions) symmetry (around $x$) in the {\stampatello short} (since $h=o(x)$) {\stampatello intervals} $[x-h,x+h]$. 
\par
For a motivation to study $I_f$, see esp. [C]. \enspace Our present arguments closely resemble [C7] ones.

\medskip

\par
Our main result is the following. 
\smallskip
\par
\noindent {\stampatello Theorem.} {\it Let } $N,h,Q\in \N$, {\it with } $Q\ll N$ {\it and } $h\to \infty$, $h=o(\sqrt N)$ {\it as } $N\to \infty$. {\it Assume } $g:\N \rightarrow \R$ {\it is independent of } $x,N,h$, {\it essentially bounded and supported in } $[1,Q]$; {\it set } $f:=g\ast \1$. {\it Then} 
$$
I_f(N,h)\EssBdd Nh. 
$$

\medskip
\par
We then get square-root cancellation for the symmetry integral of $d_k$ (generating Dirichlet series $\zeta^k$) : 
\smallskip
\par
\noindent {\stampatello Corollary.} {\it Fix } $k\ge 1$ \thinspace {\it integer}. {\it Let } $N,h\in \N$ \thinspace {\it with} \enspace $h=o(\sqrt N)$. {\it Then} ({\it the $\EssBdd-$const.depends on $k$})
$$
I_{d_k}(N,h)\EssBdd_{\!k} Nh.
$$
\par
\noindent
{\it Also, for the von Mangoldt function } $\Lambda(n):=\log p$ $\forall n=p^r$, $r\in \N$, $p$ {\it prime}, $:=0$ {\it otherwise}, 
$$
I_{\Lambda}(N,h)\EssBdd Nh.
$$

\medskip

\par
\noindent 
We'll not prove the Corollary. (A kind of \lq \lq {\stampatello Dirichlet hyperbola trick}\rq \rq \thinspace for $I_{d_k}$ gives $Q\ll N^{1-1/k}$, [C6].)

\medskip

\par
In passing, we note that an estimate of the kind in the Corollary, but for the Selberg integral of $d_k$ would give, applying [C4] bounds, the proof of the Lindel\"of Hypothesis ! (Also, we prove, with the Corollary, the almost-all symmetry of primes: for the corresponding Selberg integral, i.e. the classical Selberg integral of primes, this would give the Density Hypothesis !). 

\medskip

\par
\noindent 
The paper is organized as follows: 
\smallskip
\item{$\triangle$} in section $\S2$ we start proving the Theorem, then
\item{$\triangle$} in next section, $\S3$, we state and prove a small Lemma, in order to complete Theorem's proof;
\item{$\triangle$} we conclude with some comments, in the final section. 

\medskip

\par
\noindent
The author trusts the possibility to treat the classic Selberg integral through the approach outlined in [C]. 

\bigskip
\bigskip
\bigskip

\par
\noindent \centerline{\stampatello 2. A weak majorant principle.}
\smallskip
\par
\noindent
We start proving the Theorem, giving at next section the necessary Lemma. Here, we closely follow [C7]. 
\par
\noindent {\stampatello proof.}$\!$ Assume now on $Q\ll N$, $g:\N \rightarrow \R$, with $\supporto(g)\subset [1,Q]$, $g\ast \1:=f\EssBdd 1$, $h\to \infty$ and $h=o(N)$ : 
$$
I_f(N,h)=\sum_{x\sim N}\Big| \sum_{q\le Q,q\le x+h}g(q)\chi_q(x)\Big|^2,
$$
\par
\noindent
where we define, this time (compare [C-S], esp.), $\forall q\in \N$, 
$$
\chi_q(x)\defineq -\dashsum_{{|n-x|\le h}\atop {n\equiv 0(\!\!\bmod q)}}\sgn(n-x)=\sum_{{\ell|q}\atop {\ell>1}}{{\ell}\over q}\starsum_{j\le {{\ell}\over 2}}c^{\pm}_{j,\ell}\sin {{2\pi xj}\over {\ell}}, 
$$
\par
\noindent
and the {\stampatello Fourier coefficients are positive} 
$$
c^{\pm}_{j,q}:={4\over q}\cot {{\pi j}\over q}\sin^2 {{\pi jh}\over q}:={1\over q}F_h^{\pm}\left( {j\over q}\right)\ge 0 \quad \forall j\le {q\over 2}
$$
\par
\noindent
(better, the finite Fourier expansion has non-negative coefficients); here (see [C-S], apply Parseval identity)
$$
\starsum_{j\le q}|c^{\pm}_{j,q}|^2 \ll \sum_{0<j\le q}|c^{\pm}_{j,q}|^2 \ll \left\Vert {h\over q}\right\Vert 
\ll \min\left(1,{h\over q}\right).
$$
\par
\noindent
As usual, $\Vert \alpha \Vert:=\min_{n\in \Z}|\alpha-n|$ \thinspace is the {\it distance to the next integer}, $\forall \alpha \in \R$. 
\smallskip
\par
\noindent
The {\stampatello orthogonality of the additive characters} is the starting point (compare [C7]).
\par
Here we'll not define explicitly the {\stampatello Ramanujan coefficients}, see [C7]. 
\par
\noindent
This time, not like in [C7], we use as a \lq \lq majorant\rq \rq, for our $f(n)$, the {\stampatello divisor function}, $d(n)$. In fact, this one has, at least for $h=o(\sqrt N)$, see [C-S] (compare [C2]), square-root cancellation on average over $x\sim N$. 
\par
\noindent
Then (see [C7], proof of the Proposition: by the way, misprints occur there, missing ${1\over {\ell^2}}$ and ${1\over {\ell t}}$) 
$$
I_f(N,h)=D_f^{\pm}(N,h) 
  + 2\doublesum_{1<\ell,t\le Q}\sum_{d\le {Q\over {\ell}}}\sum_{q\le {Q\over t}}{{g(\ell d)g(tq)}\over {\ell dtq}}
     \doublestarsum_{{j\le{{\ell}\over 2} \thinspace \quad \thinspace r\le{t\over 2}}\atop {\delta:={j\over {\ell}}-{r\over t}>0}}
       F_h^{\pm}\left({j\over {\ell}}\right)F_h^{\pm}\left({r\over t}\right)
      \sum_{{x\sim N}\atop {\ldots}}\sin{{2\pi xj}\over {\ell}}\sin{{2\pi xr}\over t} = 
$$
$$
=D_f^{\pm}(N,h) 
 + \doublesum_{1<\ell,t\le Q}\sum_{d\le {Q\over {\ell}}}\sum_{q\le {Q\over t}}{{g(\ell d)g(tq)}\over {\ell dtq}}
    \doublestarsum_{{j\le{{\ell}\over 2} \thinspace \quad \thinspace r\le{t\over 2}}\atop {\delta:={j\over {\ell}}-{r\over t}>0}}
      F_h^{\pm}\left({j\over {\ell}}\right)F_h^{\pm}\left({r\over t}\right)
		\left( \sum_{{x\sim N}\atop {\ldots}}\cos 2\pi \delta x - \sum_{{x\sim N}\atop {\ldots}}\cos 2\pi \sigma x\right), 
$$
\par
\noindent
where we abbreviate 
$$
\sum_{{x\sim N}\atop {\ldots}}:=\sum_{{x\sim N}\atop {{x\ge \ell d-h}\atop {x\ge tq-h}}}, 
$$
\par
\noindent
indicating the present {\it diagonal} as 
$$
D_f^{\pm}(N,h)\defineq \sum_{1<\ell \le Q}\sum_{d_1\le {Q\over {\ell}}}\sum_{d_2\le {Q\over {\ell}}}{{g(\ell d_1)g(\ell d_2)}\over {\ell^2 d_1d_2}}
			\starsum_{j\le {{\ell}\over 2}}F_h^{\pm}\left({j\over {\ell}}\right)^2 
				\sum_{{x\sim N}\atop {{x\ge \ell d-h}\atop {x\ge tq-h}}}\sin^2 {{2\pi xj}\over {\ell}}, 
$$
\par
\noindent
with \enspace $\delta:=j/\ell - r/t$ \thinspace ($>0$, see above) and 
$$
\sigma:=\left\Vert {j\over {\ell}}+{r\over t}\right\Vert \in \left[0,{1\over 2}\right].
$$
\par
\noindent
This time we appeal to the elementary Lemma of [C7] to treat the terms with well-spaced fractions (i.e., $\delta>1/A$ or $\sigma>1/A$, with $1/A=1/NL$, $L:=\log N$, say), {\stampatello but} we have {\stampatello a minus sign} on the terms with $\sigma \le 1/A$. These are treated exactly, when $1/A=o(1/N)$, in the next small Lemma. 

\bigskip
\bigskip
\bigskip

\par
\noindent \centerline{\stampatello 3. Statement and proof of the lemma to complete the Theorem Proof.}
\smallskip
\par
We can state and show our
\smallskip
\par
\noindent {\stampatello Lemma.} {\it Let } $A,N\in \N$ {\it with } ${1\over A}=o({1\over N})$ {\it when } $N\to \infty$. {\it Let } $Q\ll N$ {\it and } $R_{\ell},R_t,F_h^{\pm}$ {\it be as above}. {\it Then}
$$
\delta:={j\over {\ell}}-{r\over t}, \enspace \sigma:=\left\Vert {j\over {\ell}}+{r\over t}\right\Vert 
\enspace \Rightarrow \enspace 
\doublesum_{1<\ell,t\le Q}R_{\ell}(f)R_t(f)
    \doublestarsum_{{j\le{{\ell}\over 2} \thinspace \quad \thinspace r\le{t\over 2}}\atop {\delta>0,\sigma \le 1/A}}
      F_h^{\pm}\left({j\over {\ell}}\right)F_h^{\pm}\left({r\over t}\right)\sum_{{x\sim N}\atop {{x\ge \ell d-h}\atop {x\ge tq-h}}}\cos 2\pi \sigma x = 0.
$$

\smallskip

\par
\noindent {\stampatello proof.}$\!$ Abbreviating $\sigma$ as above we have 
$$
\sigma \le {1\over A} \Rightarrow 0<{j\over {\ell}}+{r\over t}\le {1\over A} \enspace \hbox{\stampatello or} \enspace 0\le 1-{j\over {\ell}}-{r\over t}\le {1\over A};
$$
\par
first case gives an absurd: in particular \thinspace $0<j=o\left({{\ell}\over N}\right)=o\left({Q\over N}\right)=o(1), 0<r=o\left({t\over N}\right)=o\left({Q\over N}\right)=o(1)$. 
\par
Hence 
$$
0\le \left( {1\over 2}-{j\over {\ell}}\right)+\left( {1\over 2}-{r\over t}\right)\le {1\over A} 
\enspace \Rightarrow \enspace 
0\le {1\over 2}-{j\over {\ell}}\le {1\over A}, 0\le {1\over 2}-{r\over t}\le {1\over A},
$$
\par
whence (use $1/A=o(1/N)$, here)
$$
{{\ell}\over 2}-{{\ell}\over A}\le j \le {{\ell}\over 2}, 
\enspace 
{t\over 2}-{t\over A}\le r \le {t\over 2}, 
\enspace \Rightarrow \enspace 
j=\left[ {{\ell}\over 2}\right], r=\left[ {t\over 2}\right]
$$
\par
which, together with previous ranges for  $1-j/\ell-r/t$, give, again from ${1\over A}=o({1\over N})$, 
$$
0\le {1\over {\ell}}\left\{ {{\ell}\over 2}\right\} + {1\over t}\left\{ {t\over 2}\right\} \le {1\over A} 
\enspace \Rightarrow \enspace 
2|\ell, 2|t. 
$$
\par
Thus, 
$$
{j\over {\ell}}={r\over t}={1\over 2}
$$
\par
and the thesis, using $F_h^{\pm}\left({1\over 2}\right)=0.\enspace \square$ 

\medskip

\par
\noindent 
We may complete Theorem's proof. In fact, applying [C7] well-spaced Lemma (the conditions in $\ldots$ give a small disturbance, since the same considerations of [C7] Proposition proof apply verbatim) with the same choice \thinspace $A=NL=N\log N$, say, from present Lemma (which, together with the [C7] one, says terms with the $x-$sum of \enspace $\cos 2\pi \sigma x$ \enspace are negligible !) 
$$
I_f(N,h)\EssBdd Nh +
$$
$$
+ D_d^{\pm}(N,h) 
 + \doublesum_{1<\ell,t\le Q}\sum_{d\le {{2N+h}\over {\ell}}}\sum_{q\le {{2N+h}\over t}}{1\over {\ell dtq}}
    \doublestarsum_{{j\le{{\ell}\over 2} \thinspace \quad \thinspace r\le{t\over 2}}\atop {\delta:={j\over {\ell}}-{r\over t}>0}}
      F_h^{\pm}\left({j\over {\ell}}\right)F_h^{\pm}\left({r\over t}\right)
		\left( \sum_{{x\sim N}\atop {\ldots}}\cos 2\pi \delta x - \sum_{{x\sim N}\atop {\ldots}}\cos 2\pi \sigma x\right) 
  \EssBdd 
$$
$$
\EssBdd I_d(N,h)+Nh, 
$$
\par
\noindent
whence the thesis, since the symmetry integral of the function $d(n)$, from [C-S] (requiring $h=o(\sqrt N)$, here), is 
$$
\EssBdd Nh.\enspace \square 
$$

\bigskip
\bigskip
\bigskip

\par
\noindent \centerline{\stampatello 4. Some comments.}
\smallskip
\par
\noindent
The method outlined in [C7] and here treats more in general second moments over \lq \lq long intervals\rq \rq \thinspace of \lq \lq short interval\rq \rq \thinspace sums (or averages). Both these approaches use the non-negativity of the Fourier coefficients (say, $\widetilde{F}_h$ there and $F_h^{\pm}$ here) arising from the short interval inner average. The Selberg integral is another planet. I mean, no more free positivity there ! (No free meal...) However, there' s still some hope ! First of all, the present approach seems weaker than [C7] one, whose proofs are easier; but not a problem the estimate in our present Lemma. The real complication comes when facing changing signs of Fourier coefficients for the Selberg integral. They require another approach and, also, will set some limits on the ranges of the divisors, say $Q$. Here we don't have any troubles, since we have a REAL (though approximate !) majorant principle, i.e. we use brute force and \lq \lq bound any essentially bounded function with the divisor function\rq \rq. This is possible, since we KNOW the symmetry integral of $d(n)$, not only, but it's also fantastically small (it's a square-root saving ! However, see my paper on Parma University Journal about square-free numbers: they've even better cancellation on their symmetry integral !). 
\par
The problem we have to face for both the original and the modified Selberg integral is the (already mentioned, see [C7] final comments) problem of the Wintner majorants that even here is represented (but immediately solved).
\par
The majorating procedure is optimal (as we see now) in the case of the symmetry integral. Then, if you have a way to get non-trivial bounds from this (as Kaczorowski and Perelli do, as stated in [C]) for the Selberg integral, better to avoid a direct approach to $J_f(N,h)$.
\par
However, the modified Selberg integral has its own interest; and the further difficulty to build good Wintner majorants to get non-trivial bounds is happily faced, in view of the non-trivial results obtained !
\par
The same is of course true whenever a kind of majorant principle will be available for the Selberg integral, but (due to Fourier coefficients non-constant signs) with limitations, due to the technical analytic treatment.
\par
I expect (and have already some scratch calculations) the analytic part to be fair (not impossible, not easy); then, the \lq \lq {\it heavy}\rq \rq \thinspace part is building (non-trivial!) majorants. And next ... estimate Selberg's integral !

\vfil
\eject

\par
\noindent
\centerline{\bf References}
\smallskip
\item{\bf [C]} \thinspace Coppola, G.\thinspace - \thinspace {\sl On the symmetry of divisor sums functions in almost all short intervals} \thinspace - \thinspace Integers {\bf 4} (2004), A2, 9 pp. (electronic). $\underline{\tt MR\enspace 2005b\!:\!11153}$  
\smallskip
\item{\bf [C1]} \thinspace Coppola, G.\thinspace - \thinspace {\sl On the Correlations, Selberg integral and symmetry of sieve functions in short intervals} \thinspace - \thinspace http://arxiv.org/abs/0709.3648v3 (to appear on: Journal of Combinatorics and Number Theory)
\smallskip
\item{\bf [C2]} \thinspace Coppola, G.\thinspace - \thinspace {\sl On the Correlations, Selberg integral and symmetry of sieve functions in short intervals, II} \thinspace - \thinspace Int. J. Pure Appl. Math. {\bf 58.3}(2010), 281--298.
\smallskip
\item{\bf [C3]} \thinspace Coppola, G.\thinspace - \thinspace {\sl On the Correlations, Selberg integral and symmetry of sieve functions in short intervals, III} \thinspace - \thinspace http://arxiv.org/abs/1003.0302v1
\smallskip
\item{\bf [C4]} \thinspace Coppola, G.\thinspace - \thinspace {\sl On the Selberg integral of the $k-$divisor function and the $2k-$th moment of the Riemann zeta-function} \thinspace - \thinspace http://arxiv.org/abs/0907.5561v1 - to appear on Publ. Inst. Math., Nouv. Sér.
\smallskip
\item{\bf [C5]} \thinspace Coppola, G.\thinspace - \thinspace {\sl On the symmetry of arithmetical functions in almost all short intervals, V} \thinspace - \thinspace (electronic) http://arxiv.org/abs/0901.4738v2 
\smallskip
\item{\bf [C6]} \thinspace Coppola, G.\thinspace - \thinspace {\sl On some lower bounds of some symmetry integrals} \thinspace - \thinspace  http://arxiv.org/abs/1003.4553v2 
\smallskip
\item{\bf [C7]} \thinspace Coppola, G.\thinspace - \thinspace {\sl On the modified Selberg integral} \thinspace - \thinspace  http://arxiv.org/abs/1006.1229v1 
\smallskip
\item{\bf [C-S]} Coppola, G. and Salerno, S.\thinspace - \thinspace {\sl On the symmetry of the divisor function in almost all short intervals} \thinspace - \thinspace Acta Arith. {\bf 113} (2004), {\bf no.2}, 189--201. $\underline{\tt MR\enspace 2005a\!:\!11144}$
\smallskip
\item{\bf [D]} \thinspace Davenport, H.\thinspace - \thinspace {\sl Multiplicative Number Theory} \thinspace - \thinspace Third Edition, GTM 74, Springer, New York, 2000. $\underline{{\tt MR\enspace 2001f\!:\!11001}}$
\smallskip
\item{\bf [L]} \thinspace Linnik, Ju.V.\thinspace - \thinspace {\sl The Dispersion Method in Binary Additive Problems} \thinspace - \thinspace Translated by S. Schuur \thinspace - \thinspace American Mathematical Society, Providence, R.I. 1963. $\underline{\tt MR\enspace 29\# 5804}$
\smallskip
\item{\bf [T]} \thinspace Tenenbaum, G.\thinspace - \thinspace {\sl Introduction to Analytic and Probabilistic Number Theory} \thinspace - \thinspace Cambridge Studies in Advanced Mathematics, {\bf 46}, Cambridge University Press, 1995. $\underline{\tt MR\enspace 97e\!:\!11005b}$
\smallskip
\item{\bf [V]} \thinspace Vinogradov, I.M.\thinspace - \thinspace {\sl The Method of Trigonometrical Sums in the Theory of Numbers} - Interscience Publishers LTD, London, 1954. $\underline{{\tt MR \enspace 15,941b}}$

\medskip

\leftline{\tt Dr.Giovanni Coppola}
\leftline{\tt DIIMA - Universit\`a degli Studi di Salerno}
\leftline{\tt 84084 Fisciano (SA) - ITALY}
\leftline{\tt e-mail : gcoppola@diima.unisa.it}
\leftline{\tt e-page : www.giovannicoppola.name}

\bye